\documentclass[11pt,twoside]{article}
\title{} \author{} \date{}
\usepackage{latexsym,amssymb,times,amsmath}
\input amssym.def
\newtheorem{te}{Theorem}[section]

\newtheorem{ex}[te]{Example}

\def\bt{\hspace{-2mm}{\bf .}\hspace{2mm}}
\def\dok{\noindent{\bf Proof. }}
\def\kdok{\hfill $\Box$ \par \vspace*{2mm} }
\def\a{\alpha}
\def\b{\beta}
\def\k{\kappa}
\def\l{\lambda}
\def\m{\mu}
\def\n{\nu}
\def\o{\omega}

\def\R{{\mathbb R}}
\def\Z{{\mathbb Z}}
\def\A{{\mathcal A}}
\def\B{{\mathcal B}}
\def\C{{\mathcal C}}
\def\D{{\mathcal D}}

\def\L{{\mathcal L}}

\def\la{\langle}
\def\ra{\rangle}

\def\Jump{\mathop{{\mathrm{Jump}}}\nolimits}

\def\w{\mathop{\mbox{w}}\nolimits}
\def\cf{\mathop{\mbox{cf}}\nolimits}

\def\Init{\mathop{\mbox{Init}}\nolimits}

\begin{document}
\thispagestyle{plain}
\begin{center}
           {\large \bf WEAKLY BOOLEAN MAXIMAL CHAINS OF HOMEOMORPHIC TOPOLOGIES}
\end{center}
\begin{center}
{\bf \bf Milo\v s S.\ Kurili\'c\footnote{Department of Mathematics and Informatics, Faculty of Sciences, University of Novi Sad,
                                      Trg Dositeja Obradovi\'ca 4, 21000 Novi Sad, Serbia.
                                      e-mail: milos@dmi.uns.ac.rs}
and Bori\v sa Kuzeljevi\'c\footnote{Department of Mathematics and Informatics, Faculty of Sciences, University of Novi Sad,
                                      Trg Dositeja Obradovi\'ca 4, 21000 Novi Sad, Serbia.
                                      e-mail: borisha@dmi.uns.ac.rs }}
\end{center}
\begin{abstract}
\noindent
If $\l <\k$ are infinite cardinals,
a linear order $L$ is isomorphic to a maximal chain in $[\k ]^{\k |\k }$
(resp.\ $[\k ]^{\l |\k }$; $[\k ]^{\k |\l }$)
iff $L$ is weakly Boolean,
the weight of all initial segments of $L$ is equal to $\k$ (resp.\ $\l$, $\l ^+$) and
the weight of all final segments of $L$ is equal to $\k$ (resp.\ $\l ^+$, $\l $).
We also provide an application of this result in topology.

\vspace*{1mm}


\end{abstract}

\section{Introduction}

The purpose of this paper is to investigate maximal chains in certain subfamilies of posets of the form $(P(X),\subset)$.
For all undefined notions, please refer to Section \ref{s:notation}.
This topic has already been explored in the earlier work of, for example, Kuratowski in \cite{Kura} or Day in \cite{Day}, and the authors in \cite{K, KQ, KKG, KKP, KKB}.
The fundamental difference between the work of Kuratowski and Day on one hand, and the authors on the other hand, is that the previous work of the authors only dealt with the countable case - in a sense that we were investigating maximal chains in various subfamilies of the poset $P(\o)$.

In this paper we continue the investigation of maximal chains in subfamilies of the poset $P(\k)$ for an uncountable cardinal $\k$, in the spirit of the following already mentioned result.

\begin{te}[Kuratowski] \rm \label{T2369}
Let $\k$ be an infinite cardinal.
A linear order $L$ is isomorphic to a maximal chain in $P(\kappa )$
iff $L$ is isomorphic to $\Init (L_1 )$, for some linear order $L_1$ of size $\kappa $.
\end{te}

We will mostly deal with families where we obtain Boolean maximal chains, and this is closely related to results of the first author about maximal chains in positive families on countable sets.
Recall that a positive family on $\o$ (as defined in \cite{K}) is a family $\mathcal P\subset P(\o)$ which does not contain $\emptyset$ as an element, which is closed under taking supersets of its elements, which is closed under finite modifications of its elements, and which contains a set whose complement is countably infinite.
The central result in \cite{K} is the following (Theorem 1 in \cite{K}).

\begin{te}\label{T32011}
       Let $\mathcal P\subset P(\o)$ be a positive family. For each linear order $\langle L,< \rangle$ the following conditions are equivalent:
       \begin{enumerate}
              \item $L$ is isomorphic to some maximal chain in the poset $\langle \mathcal P\cup\{\emptyset\},\subset\rangle$;
              \item $L$ is an $\mathbb R$-embeddable complete linear order with dense jumps, and $0_L$ has no successor;
              \item $L$ is order-isomorphic to a compact nowhere dense set $K\subset [0,1]_{\mathbb R}$ such that $0\in K'$ and $1\in K$.
       \end{enumerate}
\end{te}

The structure of this paper is the following.
In Section \ref{s:notation} we introduce all the necessary notions for reading the paper.
In Section \ref{s:kappalambda} we provide characterization of order types of maximal chains in posets of the form $[\k]^{\lambda | \mu}$ where $\k,\lambda,\mu$ are infinite cardinals such that $\k=\lambda+\mu$.
In Section \ref{s:topology} we isolate a particular topology $\tau$ on $\k$ such that the set of all topologies on $\k$ homeomorphic to $\tau$ is order-isomorphic with $[\k]^{\lambda | \mu}$.
We then apply already obtained results to give a characterization of maximal chain of topologies on $\k$ homeomorphic to $\tau$.

\section{Basic notions and facts}\label{s:notation}

We use mostly standard notation.
Note that for sets $A$ and $B$, the relation $A\subset B$ means that $A$ is a subset of $B$ - it is possible that $A=B$ in this case.
For a function $f:A\to B$, and $X\subset A$ and $Y\subset B$ we write $f[X]=\{f(x):x\in X\}$ and $f^{-1}[Y]=\{x\in A: f(x)\in Y\}$.
Cardinality of a set $X$ is denoted $|X|$.
For a set $X$ and a cardinal $\k$, we define $[X]^{\k}=\{A\subset X: |A|=\k\}$.
We define $[X]^{<\k}$ and $[X]^{\le\k}$ in a similar obvious way.
For cardinals $\k,\l,\mu$ such that $\k=\l+\mu$ we also define
$$
[\k ]^{\l |\m }=\{ A\subset \k : |A|=\l \land |\k \setminus A |=\m \}.
$$

As for topological notions, we also use mostly standard notation.
For example, for a topological space $X$ and $K\subset X$, we denote by $K'$ the set of all accumulation points of the set $K$ in the space $X$.

If $\langle L,<\rangle$ is a linear order, then
the {\it end-points}, the {\it minimum} and {\it maximum} of $L$, if they
exist, will be denoted by $0$ and $1$ respectively, sometimes with appropriate subscripts.
The {\it intervals} $(x,y)_L$, $[x,y]_L$, $(-\infty , x)_L$ etc.\
are defined in the usual way. {\it Open intervals} in $L$ are the sets of the form
$(-\infty ,x)_L$, $(x,y)_L $ and $(x, \infty )_L$. Often, in order to simplify notation, we will omit the subscripts.

A pair $\langle{\mathcal A},{\mathcal B}\rangle$ is a {\it cut}
in a linear order $\langle L,<\rangle$ iff
$L={\mathcal A}\;{\buildrel {.} \over {\cup}}\; {\mathcal B}$,
${\mathcal A},{\mathcal B}\neq \emptyset$
and $a<b$, for each $a\in {\mathcal  A}$ and $b\in {\mathcal  B}$.
A cut $\langle {\mathcal A},{\mathcal B} \rangle$ is  a {\it gap} iff
neither $\max {\mathcal A}$ nor $\min {\mathcal  B}$ exist;
$\langle {\mathcal A},{\mathcal B} \rangle$ is called a {\it jump} iff
$\max {\mathcal A}$ and $\min {\mathcal  B}$ exist.

A linear order $\langle L,<\rangle $ is called:
{\it Dedekind-complete} iff there are no gaps in $\langle L,<\rangle $
(iff each non-empty subset of $L$ having an upper bound has a
{\it least upper bound (supremum)} iff each non-empty subset of $L$ having a lower bound has a
{\it greatest lower bound (infimum)});
{\it complete} iff  each subset of $L$ has a supremum
(iff each subset of $L$ has an infimum) - note that $\sup(\emptyset)=0_L$ and $\inf(\emptyset)=1_L$;
{\it dense} iff $(x,y)_L\neq \emptyset$, for each $x,y \in L$ satisfying $x<y$;
{\it ${\mathbb R}$-embeddable} iff it is isomorphic to a subset of the real line, ${\mathbb R}$;
{\it scattered} iff there is no $L_1\subset L$ such that
$L_1$ with the inherited order is a dense linear order (iff ${\mathbb Q}$ does not embed into $L$).
A linear order $\langle L,<\rangle$ is said to
{\it have dense jumps} iff there is a jump between each two elements of $L$, that is iff
for each  $a,b \in L $ satisfying $a<b$ there are $c,d \in [a,b]_L$ such that $c<d$ and $(c,d)_L=\emptyset$.
A complete linear order having dense jumps is called {\it Boolean} and
a Dedekind-complete linear order having dense jumps will be called {\it weakly Boolean}.

A subset $D$ of a linear order $\langle L,<\rangle $  is called
{\it dense} (respectively, {\it  weakly dense}) in $L$ iff for each $x,y \in L$
satisfying $x<y$ there is $z\in D$ such that $x< z < y$
(respectively, $x\leq z \leq y$). We remind the reader that, if $\langle L , < \rangle$ is a dense linear order,
then a set $D\subset L$ is weakly dense in $L$ iff it is  dense in $L$.
Let
$$
\w (L)=\min \{ |D| : D\subset L \land D \mbox{ is weakly dense in }L\}.
$$
$C\subset L$ is called {\it convex} iff for each $x,y\in C$ and each
$z\in L$ satisfying $x<z<y$ we have $z\in C$.
A subset $I$ of $L$ is said to be an {\it initial segment} of $L$ iff $L \ni x\leq y \in I$ implies $x \in I$.
It is clear that the set $\Init(L)$ of all initial segments of $L$ is linearly ordered
by the inclusion.

\section{Maximal chains in $[\k ]^{\l |\m }$}\label{s:kappalambda}

In this section we characterize the order types of maximal chains in the posets of the form
$\la [\k ]^{\l |\m }, \subset\ra$ where $\k, \l $ and $\m $ are infinite cardinals satisfying $\k = \l +\m$.
We will use the following fact.

\begin{te}[Day \cite{Day}]  \rm \label{T2366}
A linear order is isomorphic to a maximal chain in a $<\kappa$-complete atomic Boolean algebra
if and only if it is $<\kappa$-complete, has 0 and 1 and has dense jumps.
\end{te}

Note that in the following theorem, the three cases we have written down are actually all the possible cases.
This is because for infinite cardinals, if $\k=\l+\mu$, then $\l=\k$ or $\mu=\k$.

\begin{te}\bt \rm \label{T3200}
Let $\k >\l \geq \o$ be cardinals. Then
\begin{enumerate}
       \item[(a)] A linear order $L$ is isomorphic to a maximal chain in $\la [\k ]^{\k |\k }, \subset\ra$ iff $L$ is weakly Boolean and $\w ((-\infty , a])=\w ([a, \infty ))=\k$, for each $a\in L$.
       \item[(b)] A linear order $L$ is isomorphic to a maximal chain in $\la [\k ]^{\l |\k }, \subset\ra$ iff $L$ is weakly Boolean and $\w ((-\infty , a])= \l $ and $\w ([a, \infty ))=\l ^+$, for each $a\in L$.
       \item[(c)] A linear order $L$ is isomorphic to a maximal chain in $\la [\k ]^{\k |\l }, \subset\ra$ iff $L$ is weakly Boolean and $\w ((-\infty , a])= \l ^+ $ and $\w ([a, \infty ))=\l $, for each $a\in L$.
\end{enumerate}
\end{te}

\noindent
{\bf Proof of (a)}

($\Rightarrow$) Let $\L$ be a maximal chain in $[\k ]^{\k |\k }$
and $\L '$ a maximal chain in $P(\k )$ containing $\L$. By Theorem \ref{T2366} $\L '$ is
a Boolean linear order and, hence, its convex subset $\L$ is weakly Boolean.

Also we have
\begin{equation}\textstyle \label{EQ3200}
|\bigcap \L |<\k \;\; \mbox{ and } \;\; |\k \setminus \bigcup \L |<\k ,
\end{equation}
because $|\bigcap \L |=\k$ would imply that, for $\a \in \bigcap \L$, the set $\L \cup \{ (\bigcap \L )\setminus \{\a \}\}$
is a chain in $[\k ]^{\k |\k }$ extending $\L$, which is impossible by the maximality of $\L$. The assumption
$|\k \setminus \bigcup \L |=\k$ produces a similar contradiction.
For $\a \in \bigcup \L \setminus \bigcap \L$ let us define
$$\textstyle
A_\a =\bigcup \{ C\in \L : \a \not\in C\} \;\; \mbox{ and }\;\; B_\a =\bigcap \{ C\in \L : \a \in C\} .
$$
\noindent
{\bf Claim 1.}
For $\a , \b \in \bigcup \L \setminus \bigcap \L$ we have

(i) $A_\a , B_\a \in \L$

(ii) $B_\a = A_\a \cup \{ \a \}$;

(iii) $\a \neq \b \Rightarrow A_\a \neq A_\b \land B_\a \neq B_\b$.

\vspace{2mm}
\dok
(i) We have $\L =\L ^- \cup  \L ^+$, where $\L ^- = \{ C\in \L : \a \not\in C\}\neq \emptyset$
(since $\a \not\in \bigcap \L$) and $\L ^+ = \{ C\in \L : \a \in C\}\neq \emptyset$ (because $\a \in \bigcup \L$).
Since $\la \L , \subset \ra$ is a linear order, for $C_1 \in \L ^- $ and $C_2 \in \L ^+$ we have $C_1\varsubsetneq C_2$,
which implies that $A_\a = \bigcup \L ^- \subset \bigcap \L ^+ =B_\a$. These sets are comparable with each element of
$\L$ so, by the maximality of $\L$, $A_\a ,B_\a \in \L$.

(ii) Clearly we have $A_\a \cup \{ \a \} \subset B_\a$ and the proper inclusion would imply that
$\L \cup \{ A_\a \cup \{ \a \}\}$ is a chain in $[\k ]^{\k |\k }$ bigger than $\L$. Thus $A_\a \cup \{ \a \} = B_\a$.

(iii) By (ii), the equality $B_\a =B_\b$ would imply that $\a \in A_\b \setminus A_\a$ and $\b \in A_\a \setminus A_\b$, but
by (i), $A_\a , A_\b \in \L$ and, hence, $A_\a \subset A_\b$ or $A_\b \subset A_\a$.
The equality $A_\a =A_\b$ would imply that $\a ,\b \not\in A_\a $ and the sets
$B_\a = A_\a \cup \{ \a \}$ and $B_\b = A_\a \cup \{ \b \}$ would be incomparable,
which is false because $B_\a , B_\b \in \L$.
\kdok
\vspace{1mm}
\noindent
{\bf Claim 2.}
$\w ((-\infty , A]_\L )= \k$, for each $A\in \L$. In fact, if $\A =\{ A_\a : \a \in A \setminus \bigcap \L \}$,
$\B =\{ B_\a : \a \in A \setminus \bigcap \L \}$ and $\D$ is a weakly dense subset of $(-\infty , A]_\L$,
then

(i) $\A \cup \B \subset (-\infty , A]_\L$ and $|\A |=|\B |=\k$;

(ii) The set $\A$ is weakly dense in $(-\infty , A]_\L$ and, hence, $\w ((-\infty , A]_\L )\leq \k$;

(iii) $\D \cap \{ A_\a , B_\a \} \neq \emptyset$, for each $\a \in A \setminus \bigcap \L $;

(iv) $|\D |\geq \k $, thus $\w ((-\infty , A]_\L )\geq \k$.

\vspace{2mm}
\dok
(i) By Claim 1(i), $A_\a ,B_\a \in \L$ and, since $\a \in A$, $A_\a \subset B_\a \subset A$.
$|\A |=|\B | = \k$ follows from $|A|=\k$, (\ref{EQ3200}) and Claim 1(iii).

(ii) If $S,T \in (-\infty , A]_\L$ and $S\varsubsetneq T$, then, by Claim 1, for $\a \in T\setminus S$ we have
$S \subset A_\a \varsubsetneq B_\a \subset T$ so $[S,T]_\L \cap \A \neq \emptyset$.

(iii) By Claim 1(ii) we have $(A_\a , B_\a )_\L =\emptyset$. Thus, by the assumption, the set
$\D \cap [A_\a , B_\a ]_\L = \D \cap \{ A_\a , B_\a \}$ is nonempty.

(iv) By (iii) the function $f: A\setminus \bigcap \L \rightarrow \D \cap (\A \cup \B)$ given by:
$f(\a )= A_\a$, if $A_\a \in \D$; $f(\a )= B_\a$, if $A_\a \not\in \D$;
is well defined. If $f(\a )=A_\a$ for $\k$-many values of $\a$, then, by Claim 1(iii), $|\D \cap \A |=\k$.
Otherwise we have $|\D \cap \B |=\k$, thus $|\D |\geq \k $.
\kdok

Clearly, $\L ^*=\{ \k \setminus A : A\in \L \}$ is a maximal chain in $[\k ]^{\k | \k }$ anti-isomorphic to $\L$.
Hence the weight of the final parts of $\L$ is equal to the weight of the initial parts of $\L ^*$ and, applying Claim 2
to $\L ^*$ we obtain $\w ([A, \infty )_\L )= \k$, for each $A\in \L$.

\vspace{3mm}


($\Leftarrow$) Let $L$ be a weakly Boolean linear order and $\w ((-\infty , a])=\w ([a, \infty ))=\k$, for each $a\in L$.
Let $\Jump (L) =\{ \la a,b \ra \in L^2 : a<b \land (a,b)_L =\emptyset \}$.

\vspace{2mm}
\noindent
{\bf Claim 3.} Let $D=\{ b\in L : \exists a<b \;\; (a,b)_L=\emptyset \}$. Then

(i) $D$ is a weakly dense subset of $L$;

(ii) $|D|=\k$.

\vspace{2mm}
\dok
(i) Since $L$ has dense jumps, for $a,b \in L$ satisfying $a<b$ there are $c,d\in L$ such that
$a\leq c < d \leq b$ and $(c,d)_L =\emptyset$, which implies $d\in [a,b] \cap D$. Thus $D$ is weakly dense in $L$.

(ii) By the assumption we have $\w (L)=\k$ so, by (i), $|D|\geq \k$.
Let $G$ be a weakly dense subset of $L$ such that $|G|=\k$ and let  $|\Jump (L)|=\l$. For each $\la a,b \ra \in \Jump (L)$
we have $G\cap \{ a,b \} \neq \emptyset$ so the function $g: \Jump (L) \rightarrow G$ given by
$g(\la a,b \ra)=a$, if $a\in G$;  $g(\la a,b \ra)=b$, if $a\not\in G$; is well defined. Now
(since for different $\la a,b \ra ,\la a',b' \ra \in \Jump (L)$ we have $a\neq a'$ and $b\neq b'$)
if $g(\la a,b \ra)=a$, for $\l$-many $\la a,b \ra$, then $|G|\geq \l$. Otherwise,
$g(\la a,b \ra)=b$, for $\l$-many $\la a,b \ra$ and $|G|\geq \l$ again. Thus $|\Jump (L)|\leq \k$ and, since
$D=\pi _2 [\Jump (L)]$, we obtain $|D|\leq \k$.
\kdok

\vspace{2mm}
\noindent
{\bf Claim 4.} Let $\L =\{ D\cap (-\infty ,a ]_L : a \in L\}$. Then

(i) $\L \subset [D]^{\k | \k}$;

(ii) $\la \L , \varsubsetneq \ra$ is a linear order isomorphic to $L$;

(iii) $\L$ is a maximal chain in the poset $\la [D]^{\k | \k} , \subset \ra $ isomorphic to
     $\la [\k]^{\k | \k} , \subset \ra $.

\vspace{2mm}
\dok
(i) By Claim 3(i) $D\cap (-\infty ,a ]_L $ and $D\cap (a,\infty )_L $ are weakly dense subsets of $(-\infty ,a ]_L $
and $(a,\infty )_L $ respectively so, by the assumption and Claim 3(ii), they are of the size $\k$.

(ii) Clearly, the function $f:\la L, <\ra \rightarrow \la \L ,\varsubsetneq \ra$ defined by $f(a)=D\cap (-\infty ,a ]_L$
is a surjection.
If $a,b\in L$ and $a<b$ then $f(a)\subset f(b)$. Suppose that  $f(a)= f(b)$, that is $D\cap (a,b]_L=\emptyset$.
Then $\la a,b \ra$ is not a jump so there is $c\in (a,b)_L$ and, clearly, $D\cap [c,b]_L = \emptyset$,
which is false by Claim 3(i).
Thus $f$ is a strictly increasing function and, consequently, an isomorphism.

(iii) Suppose that $\L \cup \{ A \}$ is a chain in $\la [D]^{\k | \k} , \subset \ra $ and $A\not\in \L$. Let
$$
L^- =\{ a\in L : D\cap (-\infty ,a ]_L  \varsubsetneq A \} \;\mbox{ and }\;
L^+ =\{ a\in L : A\varsubsetneq D\cap (-\infty ,a ]_L  \}.
$$
Since $\bigcap _{a\in L} D\cap (-\infty ,a ]_L =\emptyset $ and $\bigcup _{a\in L} D\cap (-\infty ,a ]_L =D $ and since
$A\in [D]^{\k | \k}$, we have $L^- \neq \emptyset$ and $L^+ \neq \emptyset$, thus $\la L^- , L^+ \ra$ is a cut in $L$.

Suppose that $L^-$ has the maximum, say $a_0$. Then
\begin{equation}\textstyle \label{EQ3001}
D\cap (-\infty ,a_0 ]_L\varsubsetneq A \subset \bigcap _{a>a_0 }D\cap (-\infty ,a ]_L.
\end{equation}
If $L^+$ has the minimum, $a_1$, then
$A\varsubsetneq D\cap (-\infty ,a_1 ]_L = (D\cap (-\infty ,a_0 ]_L ) \cup \{ a_1\}$ and (\ref{EQ3001}) is impossible;
otherwise $\bigcap _{a>a_0 }D\cap (-\infty ,a ]_L = D\cap (-\infty ,a_0 ]_L$ and (\ref{EQ3001}) is impossible again.

Thus $L^-$ has no maximum and, since $L$ is a Dedekind complete linear order, $L^+$ has the minimum,
say $a_1$, and, clearly, $a_1 =\sup L^-$. Hence we have
$ \bigcup _{a<a_1} D\cap (-\infty ,a]_L= D\cap (-\infty ,a_1)_L \subset A \varsubsetneq D\cap (-\infty ,a_1 ]_L$,
which implies $a_1\in D$. But $a_1$ is not the right point of a jump. A contradiction.
\kdok


\vspace{2mm}
\noindent
{\bf Proof of (b)}

($\Rightarrow$) Let $\L$ be a maximal chain in $[\k ]^\l $. Like in (a) we show that
$\L$ is a weakly Boolean linear order and that $|\bigcap \L |<\l$ and $|\bigcup \L |\geq \l ^+$, which implies
$\cf (\L )\geq \l ^+$.
Since a chain of sets having the cofinality $\geq \l ^{++}$ must have a member of size $\l ^+$, we have
$\cf (\L )= \l ^+$, which implies $|\bigcup \L |= \l ^+$. Thus
\begin{equation}\textstyle \label{EQ3202}
|\bigcap \L |<\l \;\; \mbox{ and } \;\; \cf (\L )=|\bigcup \L |= \l ^+ .
\end{equation}
\noindent
{\bf Claim 5.}
$\w ((-\infty , A]_\L )= \l$, for each $A\in \L$.

\vspace{2mm}
\dok
Let $\la A_\a : \a < \l ^+ \ra$ be a strictly increasing cofinal sequence in $[A, \infty )_\L$. Then
$|A_\l \setminus A|=\l$ and, hence, $A \in [A_\l ]^{\l | \l }$. Since $(-\infty , A]_\L$ is a chain in
$[A_\l ]^{\l | \l }$, there is a maximal chain $\L '$ in $[A_\l ]^{\l | \l }$ such that $(-\infty , A]_\L \subset \L '$ and
we prove that
\begin{equation}\textstyle \label{EQ3203}
(-\infty , A]_\L =(-\infty , A]_{\L '}.
\end{equation}
The inclusion ``$\subset$" is trivial. If $C\in \L '$ and $C\subset A$, then,
since $(-\infty , A]_\L \subset \L '$, $C$ is comparable with all elements of $(-\infty , A]_\L$ and, clearly, with
all elements of $[A, \infty )_\L$. Thus $\L \cup \{ C \}$ is chain in $[\k ]^\l$, and, by the maximality of $\L$,
$C\in \L$. The inclusion ``$\supset$" in (\ref{EQ3203}) is proved. Applying (a) of this theorem to
the maximal chain $\L '$ in $[A_\l ]^{\l | \l }$ we obtain $\w ((-\infty , A]_{\L '})= \l$ and, by (\ref{EQ3203}),
$\w ((-\infty , A]_\L )= \l$.
\kdok

\vspace{2mm}
\noindent
{\bf Claim 6.}
$\w ([A, \infty )_\L )= \l ^+$, for each $A\in \L$.

\vspace{2mm}
\dok
By (\ref{EQ3202}) we have $\cf ([A, \infty )_\L )=\cf (\L )=\l ^+$ and, hence, $\w ([A, \infty )_\L )\geq \l ^+$.
Let $\la A_\a : \a < \l ^+ \ra$ be a strictly increasing cofinal sequence in $\L$. By Claim 5, for each $\a < \l ^+$
there is a weakly dense subset $D_\a$ of $(-\infty , A_\a ]_\L$, such that $|D_\a |=\l$. Now
$D=\bigcup _{\a <\l ^+}D_\a$ is a weakly dense subset of $\L$ and $|D|\leq \l ^+$, which implies
$\w ([A, \infty )_\L )\leq \l ^+$.
\kdok


\vspace{2mm}
($\Leftarrow$)  Let $L$ be a weakly Boolean linear order and let $\w ((-\infty , a])=\l $ and
$\w ([a, \infty ))=\l ^+ $, for each $a\in L$.
Let $\Jump (L) =\{ \la a,b \ra \in L^2 : a<b \land (a,b)_L =\emptyset \}$.

\vspace{2mm}
\noindent
{\bf Claim 7.} Let $D=\{ b\in L : \exists a<b \;\; (a,b)_L=\emptyset \}$. Then

(i) $D$ is a weakly dense subset of $L$;

(ii) $|D|=\l ^+$.

\vspace{2mm}
\dok
(i) See the proof of Claim 3.

(ii) By the assumption we have $\w (L)=\l ^+$ so, by (i), $|D|\geq \l ^+$.
Let $G$ be a weakly dense subset of $L$ such that $|G|=\l ^+$.
As in the proof of Claim 3(ii) we show that
$|\Jump (L)|\leq |G|$ and $|D|\leq \l ^+$.
\kdok

\vspace{2mm}
\noindent
{\bf Claim 8.} Let $\L =\{ D\cap (-\infty ,a ]_L : a \in L\}$. Then

(i) $\L \subset [D]^{\l | \l ^+}$;

(ii) $\la \L , \varsubsetneq \ra$ is a linear order isomorphic to $L$;

(iii) $\L$ is a maximal chain in the poset $\la [D]^{\l | \l ^+} , \subset \ra $.

\vspace{2mm}
\dok
(i) By Claim 7(i) $D\cap (-\infty ,a ]_L $ is a weakly dense subset of $(-\infty ,a ]_L $
so, by the assumption, $|D\cap (-\infty ,a ]_L| \geq \l$.
Let $G$ be a weakly dense subset of $(-\infty ,a ]_L$ such that $|G|=\l $.
As in the proof of Claim 3(ii) we show that $|D\cap (-\infty ,a ]_L| \leq |\Jump (L)\cap (-\infty ,a ]_L|\leq |G|=\l$.
Thus $|D\cap (-\infty ,a ]_L|=\l$ and, since $|D|=\l ^+$, $|D\setminus (D\cap (-\infty ,a ]_L )|=\l ^+$.

(ii) See the proof of Claim 4(ii).

(iii) Suppose that $\L \cup \{ A \}$ is a chain in $\la [D]^{\l | \l ^+} , \subset \ra $ and $A\not\in \L$. Let
$$
L^- =\{ a\in L : D\cap (-\infty ,a ]_L  \varsubsetneq A \} \;\mbox{ and }\;
L^+ =\{ a\in L : A\varsubsetneq D\cap (-\infty ,a ]_L  \}.
$$
Since $\bigcap _{a\in L} D\cap (-\infty ,a ]_L =\emptyset $ and $|\bigcup _{a\in L} D\cap (-\infty ,a ]_L |=|D|=\l ^+ $ and since
$A\in [D]^{\l }$, we have $L^- \neq \emptyset$ and $L^+ \neq \emptyset$, thus $\la L^- , L^+ \ra$ is a cut in $L$.
The rest of the proof is identical to the corresponding part of the proof of Claim 4(iii).
\kdok

\vspace{2mm}
By Claim 7(ii), there is a bijection $f: D\rightarrow \l ^+ \subset \k$. Let the function
$F:P(D)\rightarrow P(\l ^+ )$ be defined by $F(S)=f[S]$ and let
$$\L _1 = F[\L ] = \{ f[D\cap (-\infty ,a ]_L] : a \in L \}.$$
\noindent
{\bf Claim 9.}
(i) $F[[D]^\l ]= [\l ^+ ]^\l \subset [\k ]^\l$;

(ii) $F\mid [D]^\l $ is an isomorphism of the posets $\la [D]^\l  , \subset \ra $ and $\la [\l ^+ ]^\l  , \subset \ra $;

(iii) $\L _1 $ is a maximal chain in $\la [\l ^+ ]^\l  , \subset \ra $ isomorphic to $\L$ and, hence, to $L$;

(iv) $\L _1 $ is a maximal chain in the poset $\la [\k ]^\l  , \subset \ra $.

\vspace{2mm}
\dok
Assertions (i) and (ii) are evident and (iii) follows from (ii) and Claim 8.

(iv) Suppose that $A\in [\k ]^\l \setminus \L _1$ and that $\L _1 \cup \{ A \}$ is a chain. Let
$$
L^- =\{ a\in L : f[D\cap (-\infty ,a ]_L ] \varsubsetneq A \}, \;
L^+ =\{ a\in L : A\varsubsetneq f[D\cap (-\infty ,a ]_L ] \}.
$$
Using (iii), as in the proof of (\ref{EQ3202}) we obtain $|\bigcap \L _1|<\l $ and  $|\bigcup \L _1|=\l ^+ $
and, hence, $L^- , L^+ \neq \emptyset$. But this implies that $A\in [\l ^+] ^\l$, which is impossible, by the maximality of
$\L _1$ in $[\l ^+] ^\l$.
\kdok

\vspace{2mm}
\noindent
{\bf Proof of (c)}

The mapping $f: [\k ]^{\l | \k }\rightarrow [\k ]^{\k | \l }$ defined by $f(A)=\k \setminus A$
is an isomorphism of the posets $\la [\k ]^{\l | \k }, \subset \ra$ and $\la [\k ]^{\k | \l }, \supset \ra$ so
(c) follows from (b).
\kdok

In particular, for maximal chains in the poset $\la [\o ]^{\o | \o}, \subset \ra$ we have

\begin{te}\bt \rm \label{T3202}
For each linear order $L$ the following conditions are equivalent:
\begin{enumerate}
       \item[(a)] $L$ is isomorphic to a maximal chain in the poset $\la [\o ]^{\o | \o}, \subset \ra$;
       \item[(b)] $L$ is $\R$-embeddable, weakly Boolean and without end-points;
       \item[(c)] $L \cong K\setminus \{ 0,1 \}$, for some compact nowhere dense set $K\subset [0,1]_\R$ satisfying $0,1\in K'$.
\end{enumerate}
\end{te}

\dok
The equivalence (a) $\Leftrightarrow$ (b) follows from Theorem \ref{T3201}(a). Namely, first,
the initial and final segments of $L$ have weight $\aleph _0$ iff $L$ is of weight $\aleph _0$ and without end points.
Second, it is well known that $L$ has weight $\aleph _0$ iff it is $\R$-embeddable.

(a) $\Leftrightarrow$ (c) From Theorem \ref{T32011} it follows that
a linear order $L$ is isomorphic to a maximal chain in the poset $\la [\o ]^\o , \subset \ra$ iff
$L \cong K\setminus \{ 0 \}$, for some compact nowhere dense set $K\subset [0,1]_\R$ such that $0\in K'$ and $1\in K$.
Now, $\L$ is a maximal chain in the poset $[\o ]^{\o | \o}$ iff for some (each) $A\in \L$ both
\begin{enumerate}
       \item[-] the interval $(-\infty ,A]_\L$ is a maximal chain in $[A]^\o \cong [\o ]^\o $ and, hence, isomorphic to $K_0 \setminus \{0\}$ for a compact nowhere dense $K_0\subset [0,1]_\R$ such that $0\in K_0 '$ and $1\in K_0$, and
       \item[-] the interval $[A,\infty )_\L$ is anti-isomorphic (via complementation)
       to a maximal chain in $[\o \setminus A]^\o \cong [\o ]^\o $ and, hence, anti-isomorphic to
       $K_1 \setminus \{0\}$ for a compact nowhere dense $K_1\subset [0,1]_\R$ such that $0\in K_1 '$ and $1\in K_1$;
\end{enumerate}
iff $\L$ is isomorphic to the sum $(K_0 \setminus \{0\}) + (K_1 \setminus \{0\} )^*$, which is, clearly, isomorphic to a linear order of the form $K\setminus \{ 0,1 \}$ described in (c).
\kdok

\begin{ex}\bt\rm \label{EX3200}
Simple maximal chains in $[\k ]^{\l |\n }$.
By Theorem \ref{T3202}, in $[\o ]^{\o |\o }$ there are maximal chains isomorphic to:
$\C \setminus \{ 0,1 \}$, where $\C \subset [0,1]_\R$ is the Cantor ternary set;
$\a ^* +\b$, where $\a$ and $\b$ are countable limit ordinals;
in particular, to $\Z \cong \o ^* + \o$.
By Theorem \ref{T3200}, if $\k >\l \geq \o$, then the linear order
$\k ^* + \k$ (resp.\ $\l ^* + \l ^+$, $(\l ^+)^* + \l $) is isomorphic to a maximal chain in
$[\k ]^{\k |\k }$ (resp.\ $[\k ]^{\l |\k }$, $[\k ]^{\k |\l }$).
The sum $\sum _{\k ^* + \k }\C$ is isomorphic to a maximal chain in
$[\k ]^{\k |\k }$, etc.
\end{ex}

Regarding the size of its elements, the poset $P(\k )$ naturally splits
into convex classes. For example
\begin{equation}\label{EQ3204}\textstyle
P(\o _1 )= [\o _1]^{<\o } \cup [\o _1]^{\o } \cup [\o _1]^{\o _1 |\o _1} \cup [\o _1]^{\o _1 |\o } \cup
[\o _1]^{\o _1 |<\o} .
\end{equation}
In the sequel we give some comments concerning the intersections of maximal chains with these classes.

\begin{ex}\bt\rm \label{EX3201}
Maximal chains intersecting all classes from (\ref{EQ3204}).
If $L=\o _1 +\o _1^*$, then, by Kuratowski's theorem, $\Init (L)\cong \o _1 +1 + \o _1 ^*$ is a maximal chain in $P(L)$.
If $f:L\rightarrow \o _1$ is a bijection, then $\{ f[I]: I\in \Init (L)\}$ is a maximal chain in $P(\o _1)$
intersecting all the classes listed in (\ref{EQ3204}). Note that some of these intersections are {\it not}
maximal chains in the corresponding classes.
\end{ex}
\begin{ex}\bt\rm \label{EX3202}
Maximal chains in $P(\k )\setminus \{ \emptyset , \k \}$ contained in one class.
If $L=\o _1 ^* +\o _1$, then, by Theorem \ref{T2369}, $\Init (L)\cong (\o _1 +1 )^* + \o _1 +1$ is a maximal chain in $P(L)$.
If $f:L\rightarrow \o _1$ is a bijection, then $\{ f[I]: I\in \Init (L) \setminus \{ \emptyset , L \}\}$
is a maximal chain in $P(\o _1)\setminus \{ \emptyset , \o _1 \}$ contained in $[\o _1]^{\o _1 |\o _1}$
and isomorphic to $\o _1 ^* +\o _1$.

If $L=\o  ^* +\o _1$, then $\Init (L)\cong (\o  +1 )^* + \o _1 +1$ is a maximal chain in $P(L)$ again.
If $f:L\rightarrow \o _1$ is a bijection, then $\{ f[I]: I\in \Init (L) \setminus \{ \emptyset , L \}\}$
is a maximal chain in $P(\o _1)\setminus \{ \emptyset , \o _1 \}$ contained in $[\o _1]^{\o |\o _1}$
and isomorphic to $\o  ^* +\o _1$.
\end{ex}

\section{An application to topology}\label{s:topology}

In this section we apply Theorem \ref{T3200}, and describe some classes of homeomorphic topologies on a cardinal $\k$. Recall that $\operatorname{Top}_{\k}$ denotes the set of all topologies on $\k$, that for two topologies $\tau$ and $\sigma$ we write $\tau\cong\sigma$ to denote that they are homeomorphic, and that for $\tau\in\operatorname{Top}_{\k}$ we also denote $$[\tau]_{\cong}=\{\sigma\in\operatorname{Top}_{\k}:\ \tau\cong\sigma\}.$$

\begin{te}\bt \rm \label{T3201}
Let $\k,\l,\mu$ be infinite cardinals such that $\k=\l+\mu$.
Then there is a topology $\tau$ on $\k$ such that the following two conditions hold:
\begin{enumerate}
    \item[(a)] $\left([\tau]_{\cong},\subset\right)\cong\left([\k]^{\lambda | \mu},\subset\right)$ (note that this is an isomorphism of posets);
    \item[(b)] For any linear order $L$:
    \begin{enumerate}
        \item[(i)] if $\k=\l=\mu$, then there is a maximal chain $\mathcal L\in [\tau]_{\cong}$ such that $\mathcal L\cong L$ if and only if $L$ is weakly Boolean linear order and $\w((\cdot,a])=\w([a,\cdot))=\k$ for each $a\in L$
        \item[(ii)] if $\l<\k$, then there is a maximal chain $\mathcal L\in [\tau]_{\cong}$ such that $\mathcal L\cong L$ if and only if $L$ is weakly Boolean linear order and both $\w((\cdot,a])=\l$ and $\w([a,\cdot))=\l^+$ hold for each $a\in L$
        \item[(iii)] if $\mu<\k$, then there is a maximal chain $\mathcal L\in [\tau]_{\cong}$ such that $\mathcal L\cong L$ if and only if $L$ is weakly Boolean linear order and both $\w((\cdot,a])=\mu^+$ and $\w([a,\cdot))=\mu$ hold for each $a\in L$
    \end{enumerate}
\end{enumerate}
\end{te}

\dok
Let us define a topology $\tau$ on $\k$ in the following way.
Take any $S\subset \k$ such that $ |  S |  =\lambda$ and $ |  \k\setminus S | =\mu$.
Now define $\tau=P(S)\cup \{\k\}$.
Note that $\tau\in\operatorname{Top}_{\k}$ since clearly $\tau$ is defined on $\k$, both $\emptyset$ and $\k$ belong to $\tau$ as elements, since arbitrary union of elements in $\tau$ is still in $\tau$, and since intersection of arbitrary elements of $\tau$ is again in $\tau$.
We will show that $\tau$ satisfies both conclusions from the statement of the theorem.

Proof of (a). To see this first note that if $\sigma$ is in $\operatorname{Top}_{\k}$ and $\sigma\cong\tau$, this means that there is a bijection $f:\k\to\k$ such that $\sigma=\{f[O]:O\in \tau\}$. On the other hand, for any bijection from $\k$ to $\k$ we have that $\{f[O]:O\in\tau\}=P(f[S])\cup\{\k\}$, which is clearly a topology on $\k$ homeomorphic to $\tau$.
Hence, if for a topology $\sigma\in \operatorname{Top}_{\k}$ and $f:\k\to\k$ we define $f^*(\sigma)=\{f[O]:O\in\sigma\}$, then we have $$[\tau]_{\cong}=\{f^*(\tau):f\in \operatorname{Sym}(\k)\}.$$
Let us denote $A=\{f^*(\tau):f\in \operatorname{Sym}(\k)\}$. We will prove that $(A,\subset)$ is isomorphic to $[\k]^{\l | \mu}$ as partially ordered sets.
Define $\Phi:A\to [\k]^{\l | \mu}$ as follows: for $f\in \operatorname{Sym}(\k)$ let $\Phi(f^*(\tau))=f[S]$.
Since each $f\in \operatorname{Sym}(\k)$ is a bijection from $\k$ onto $\k$ and $S\in [\k]^{\l | \mu}$, it is clear that $f[S]\in [\k]^{\l | \mu}$ as well.
Hence, $\Phi$ is a well defined map.
Next, we show that $\Phi$ is a bijection.
Suppose that $f^*(\tau)\neq g^*(\tau)$ for some $f,g\in \operatorname{Sym}(\k)$.
By the observation from the beginning of the proof, then we have $$P(f[S])\cup\{\k\}=f^*(\tau)\neq g^*(\tau)=P(g[S])\cup\{\k\}.$$
Thus, $f[S]\neq g[S]$, i.e. $\Phi(f^*(\tau))\neq \Phi(g^*(\tau))$ and so $\Phi$ is 1-1.
Now take any $T\in [\k]^{\l | \mu}$.
Since $ |  S | = |  T | =\l$ there is a bijection $f_1:S\to T$.
Since $ |  \k\setminus S | = |  \k\setminus T | =\mu$ there is a bijection $f_2:\k\setminus S\to \k\setminus T$.
Now it is clear that $f=f_1\cup f_2$ is bijection from $\k$ to $\k$ such that $f[S]=T$, i.e. $$\Phi(f^*(\tau))=f[S]=T,$$
so $\Phi$ is onto as well.
We still have to show that $\Phi$ preserves the order.
Let $f^*(\tau)\subset g^*(\tau)$ for $f,g\in \operatorname{Sym}(\k)$.
Then $$\{f[O]:O\in \tau\}=f^*(\tau)\subset g^*(\tau)=\{g[O]:O\in\tau\}.$$
In particular, since $S\in \tau$ there is some $O\in \tau$ such that $g[O]=f[S]$.
By the definition of $\tau$, there are two options, either $O=\k$ or $O\subset S$.
Note that $O=\k$ is not possible since $g[\k]=\k\neq f[S]$.
Hence, $O\subset S$.
Then $$\Phi(f^*(\tau))=f[S]=g[O]\subset g[S]=\Phi(g^*(\tau)),$$
as required.

Proof of (b).
Follows directly from (a) together with Theorem \ref{T3200} and the observation that if $\k=\l+\mu$, then $\l$ or $\mu$ has to be equal to $\k$.
\kdok
\paragraph{Acknowledgement.}
This research was supported by the Science Fund of the Republic of Serbia,
Program IDEAS, Grant No.\ 7750027:
{\it Set-theoretic, model-theoretic and Ramsey-theoretic
phenomena in mathematical structures: similarity and diversity}--SMART.


\footnotesize

\begin{thebibliography}{22}
\bibitem{Day}
      G.\ W.\ Day,
      Maximal chains in atomic Boolean algebras,
      Fund.\ Math.\ 67 (1970) 293--296.
\bibitem{Kura}
      K.\ Kuratowski,
      Sur la notion de l'ordre dans la th\'eorie des ensembles,
      Fund.\ Math.\ 2 (1921) 161--171.
\bibitem{K}
      M.\ S.\ Kurili\'c,
      Maximal chains in positive subfamilies of $P(\omega )$, 
      Order 29,1 (2012) 119--129.
\bibitem{KQ}
       M.\ S.\ Kurilić,
       Maximal chains of copies of the rational line.
       Order 30 (2013) 737--748.
\bibitem{KKP}
       M.\ S.\ Kurilić, B.\ Kuzeljević,
       Maximal chains of isomorphic suborders of countable ultrahomogeneous partial orders,
       Order 32 (2015) 83--99.
\bibitem{KKG}
       M.\ S.\ Kurilić, B.\ Kuzeljević,
       Maximal chains of isomorphic subgraphs of countable ultrahomogeneous graphs,
       Adv.\ Math.\ 264 (2014) 762--775.
\bibitem{KKB}
       M.\ S.\ Kurilić, B. Kuzeljević,
       Positive families and Boolean chains of copies of ultrahomogeneous structures,
       C.\ R.\ Math.\ Acad.\ Sci.\ Paris 358 (2020) 791--796.
\end{thebibliography}
\end{document}